\documentclass[11pt]{article}
\usepackage{overcite}
\usepackage{graphicx}
\usepackage{mathrsfs}
\usepackage{latexsym}
\usepackage{psfrag}
\usepackage{graphicx,subfigure}
\usepackage{fancyhdr,graphicx}
\usepackage{multicol}
\usepackage{cite}
\usepackage{amsmath,amssymb}
\usepackage{flafter}
\usepackage{fancyhdr}
\usepackage{subfigure}
\usepackage{indentfirst,latexsym,bm}

\usepackage{xcolor}  

\usepackage{algorithm} 
\usepackage{algorithmic} 

\newtheorem{theorem}{Theorem}[section]
\newtheorem{lemma}[theorem]{Lemma}
\newtheorem{corollary}[theorem]{Corollary}

\newtheorem{problem}{Problem}
\newtheorem{example}[theorem]{Example}

\newenvironment{pf}{
\par
\noindent {\bf Proof.}\rm}
{\mbox{}\hfill\rule{0.5em}{0.8em} \par \bigskip}

\begin{document}

\title{\bf On $t$-relaxed coloring of complete multi-partite graphs
\thanks{Supported by NSFC No. 11771080.}}

\author{Jun Lan$^{a}$  and  Wensong Lin$^{a,}$\footnote{Corresponding author. E-mail address: wslin@seu.edu.cn}\\
{\small $^{a}$ School of Mathematics, Southeast University, Nanjing 210096, P.R. China}}

\date{}
\maketitle

\vspace*{-1cm} \setlength\baselineskip{7mm}
\bigskip

\begin{abstract}
Let $G$ be a graph and $t$ a nonnegative integer. Suppose $f$ is a mapping from the vertex set of $G$ to $\{1,2,\dots, k\}$. If, for any vertex $u$ of $G$, the number of neighbors $v$ of $u$ with $f(v)=f(u)$ is less than or equal to $t$, then $f$ is called a $t$-relaxed $k$-coloring of $G$. And $G$ is said to be $(k,t)$-colorable. The $t$-relaxed chromatic number of $G$, denote by $\chi_t(G)$, is defined as the minimum integer $k$ such that $G$ is $(k,t)$-colorable.  A set $S$ of vertices in $G$ is $t$-sparse if $S$ induces a graph with a maximum degree of at most $t$. Thus $G$ is $(k,t)$-colorable if and only if the vertex set of $G$ can be partitioned into $k$ $t$-sparse sets.  It was proved by Belmonte, Lampis and Mitsou (2017) that the problem of deciding if a complete multi-partite graph is $(k,t)$-colorable is NP-complete. In this paper, we first give tight lower and up bounds for the $t$-relaxed chromatic number of complete multi-partite graphs. And then we design an algorithm to compute maximum $t$-sparse sets of complete multi-partite graphs running in $O((t+1)^2)$ time. Applying this algorithm, we show that the greedy algorithm for $\chi_t(G)$ is $2$-approximate and runs in $O(tn)$ time steps (where $n$ is the vertex number of $G$). In particular, we prove that for $t\in \{1,2,3,4,5,6\}$, the greedy algorithm produces an optimal $t$-relaxed coloring of a complete multi-partite graph. While, for $t\ge 7$, examples are given to illustrate that the greedy strategy does not always construct an optimal $t$-relaxed coloring. \\



\end{abstract}

\section{Introduction}

\noindent Suppose there are many transmitters in an area. We need to assign a channel to each transmitter through which it can transmit signals. To avoid interference, ``close" transmitters are required to receive different channels. The problem is to assign channels to transmitters so that the number of channels assigned is as small as possible. This channel assignment problem can be modeled as the graph coloring problem. We use vertices of a graph to denote the transmitters and if two transmitters are ``close" to each other then the corresponding two vertices are adjacent. The graph so obtained is referred as the interference graph of the channel assignment problem. To minimize the number of channels assigned to transmitters is equivalent to determine the chromatic number of the interference graph.

In practice, sometimes, low level of interference is acceptable. Thus the requirement that ``close" transmitters must receive different channels can be relaxed in some extend. This leads to the concept of relaxed (or improper, or defective) coloring of a graph. Please refer to [\citeonline{ABGHMM12,CCW86,GMS16}].

Let $G$ be a graph and $t$ a nonnegative integer. Suppose $f$ is a mapping from the vertex set of $G$ to $\{1,2,\dots, k\}$. If, for any vertex $u$ of $G$, the number of neighbors $v$ of $u$ with $f(v)=f(u)$ is less than or equal to $t$, then $f$ is called a {\em $t$-relaxed $k$-coloring} of $G$. And $G$ is said to be {\em $(k,t)$-colorable}. The {\em $t$-relaxed chromatic number} of $G$, denote by $\chi_t(G)$, is defined as the minimum integer $k$ such that $G$ is $(k,t)$-colorable. It is clear that $\chi_0(G)$ is the same as $\chi(G)$, the chromatic number of $G$.

Let $S$ be a subset of $V(G)$. If the induced subgraph $G[S]$ has the maximum degree at most $t$ then we call $S$ a {\em $t$-sparse set} of $G$. Let $\beta_t(G)$ denote the order of the largest $t$-sparse set of $G$. It is obvious that in a $t$-relaxed coloring of $G$ each color class is a $t$-sparse set. Thus $\chi_t(G)\ge |V(G)|/ \beta_t(G)$ for any graph $G$. For a graph $G$, let $\Delta(G)$ denote the maximum degree of $G$. Let $\omega(G)$ denote the {\em clique number}  of $G$ (the number of vertices in a maximum complete subgraph of $G$). The Brooks type result about the $t$-relaxed chromatic number can be derived from a result proved by Lov\'{a}sz in [\citeonline{Lova66}].

\begin{theorem} \label{Lova}  [\citeonline{Lova66}] \
Let $G$ be a graph and $k_1,k_2,\dots,k_m$ nonnegative integers such that $\Delta(G)\le \sum\limits_{i=1}^{m} (k_i+1)-1$, then there is a partition of $V(G)$ into $m$ classes $V_1,V_2,\dots,V_m$ such that $\Delta(G[V_i])\le k_i$ for $i=1,2,\dots,m$, where $G[V_i]$ is the subgraph of $G$ induced by the vertex set $V_i$.
\end{theorem}

By letting $k_i=t$ for $i=1,2,\dots,m$, we see that if $m\ge \lceil \frac{\Delta(G)+1}{t+1} \rceil$ then $G$ is $(m,t)$-colorable. This yields the following upper bound for $t$-relaxed chromatic number of a graph.

\begin{corollary} \label{ub}
For any graph $G$, $\chi_t(G)\le \lceil \frac{\Delta(G)+1}{t+1} \rceil$.
\end{corollary}

There is also a relation between $\chi_t(G)$ and $\chi(G)$.

\begin{theorem}\label{lb}
For any graph $G$, $\chi_t(G)\ge \lceil \frac{\chi(G)}{t+1} \rceil$.
\end{theorem}

\begin{pf}
Let $k=\chi_t(G)$. Suppose $f$ is a $t$-relaxed coloring of $G$ that partitions $V(G)$ into $k$ color classes $V_1,V_2,\dots,V_k$. Then, for each $i=1,2,\dots,k$, the maximum degree of the subgraph $G[V_i]$ (the subgraph induced by $V_i$) is at most $t$. By the Brook's Theorem, $\chi(G[V_i])\le t+1$. That is, each $V_i$ can be decomposed into at most $t+1$ independent sets. It follows that $\chi(G)\le k(t+1)=\chi_t(G) (t+1)$, implying $\chi_t(G)\ge \frac{\chi(G)}{t+1}$. Since $\chi_t(G)$ is integral, the theorem follows. 
\end{pf}

Denote by $K_n$ the complete graph on $n$ vertices. From Corollary \ref{ub} and Theorem \ref{lb}, we have $\chi_t(K_n)=\lceil \frac{n}{t+1} \rceil$. Thus both the upper bound for $\chi_t(G)$ in Corollary \ref{ub} and the lower bound in Theorem \ref{lb} are attainable.

For any two fixed integers $k\ge 2$ and $t\ge 1$, it was proved in [\citeonline{CGJ97}] that the problem of deciding if a graph is $(k,t)$-colorable is NP-complete. All planar graphs are proved to be $(3,2)$-colorable in [\citeonline{CCW86}]. And in [\citeonline{CGJ97}], the authors proved the NP-completeness of the $(2,t)$-colorable problem ($t\ge 1$) and $(3,1)$-colorable problem for planar graphs. Havet, Kang, and  Sereni in [\citeonline{HKS09}]  showed that the problem of deciding if a unit disk graph is $(k,t)$-colorable is NP-complete for any two fixed integers $k\ge 3$ and $t\ge 0$. A graph is an {\em interval graph} if its vertices are closed intervals of $\mathbb{R}$ and two vertices are adjacent if the corresponding two intervals intersect. For any two fixed nonnegative integers $k$ and $t$, it was proved in [\citeonline{HKS09}] that the problem of deciding if an interval graph is $(k,t)$-colorable is solvable in polynomial time. A {\em unit interval graph} is an interval graph with every interval having unit length. It was proved in [\citeonline{HKS09}] that, for any unit interval graph $G$, $\lceil \frac{\omega(G)}{t+1} \rceil \le \chi_t(G)\le \lceil \frac{\omega(G)}{t+1} \rceil +1$. However, as indicated by the authors in [\citeonline{HKS09}], even for unit interval graphs, it is not known, for fixed integer $t\ge 1$, whether there is a polynomial time algorithm to determine its $t$-relaxed chromatic number.
For any two positive integers $n$ and $r$ with $n\ge r+1$, denote by $P_n^r$ the $r$th power of a path on $n$ vertices. Let $t$ be any positive integer, the $t$-relaxed chromatic number of $P_n^r$ was determined in [\citeonline{LL2018}]. 

It was proved by Belmonte etc. in [\citeonline{BLM2017}] that deciding if a complete multi-partite graph is $(k,t)$-colorable is a NP-complete problem. Based on dynamic programming method, they showed that there is an algorithm which decides if a cograph admits a $(k,t)$-coloring in time $O^* (k^{O(t^4)})$. Note that complete multi-partite graphs are cographs. Thus if $t$ is bounded then the $(k,t)$-colorable problem in cographs is polynomial-time solvable.

This paper focuses on $t$-relaxed colorings of complete multi-partite graphs.
Tight lower and up bounds for the $t$-relaxed chromatic number of complete multi-partite graphs are given in next section.  In Section 3, we first design an algorithm to compute the maximum $t$-sparse set of a complete multi-partite graph running in $O((t+1)^2)$ time. Applying this algorithm, we then describe the greedy algorithm that compute a $t$-relaxed coloring of a complete multi-partite graph in $O(tn)$ time steps (where $n$ is the vertex number of a graph). The approximation ratio of the greedy algorithm is proved to be at most $2$.  In Section 4 and 5, we show that the greedy algorithm actually outputs an optimal $t$-relaxed coloring of a complete multi-partite graph for $t=1,3$, respectively. The similar arguments in Section 4 and 5 work for the case when $1\le  t \le 6$.  However, if $t\ge 7$ they do not work any more. In Section 6, for $t\ge 7$, examples are given to illustrate that the greedy strategy does not always construct an optimal $t$-relaxed coloring.

\section{Lower and upper bounds}

Let $s\ge 2$ be an integer. Let $n_1,n_2,\dots,n_s$ be $s$ positive integers. For convenience, we always assume $n_1\ge n_2\ge \cdots \ge n_s$. Denote by  $K(n_1,n_2,\dots,n_s)$ the complete $s$-partite graph with the $s$ parts having $n_1,n_2, \dots, n_s$ vertices, respectively. For $i=1,2,\dots,s$, let $V_i$ be the set of vertices in the $i$th part of $K(n_1,n_2,\dots,n_s)$. It is obvious that $\chi(K(n_1,n_2,\dots,n_s))=s=\omega(K(n_1,n_2,\dots,n_s))$.

In this paper, we always assume that $t$ is a fixed positive integer.  Suppose $f$ is a $t$-relaxed coloring of $K(n_1,n_2,\dots,n_s)$. For any color $i$ used by $f$, denote by $f^{-1}(i)$ the set of vertices that are assigned the color $i$. We say the color $i$ appears on some part $V_j$ if there is at least one vertex in $V_j$ assigned the color $i$. Let $S$ be a subset of $V(G)$, denote by $f(S)$ the set of colors that are assigned to vertices in $S$. Let $||f||$ denote the order of $f(V(G))$, that is, the number of distinct colors used by $f$.

\begin{lemma}\label{mvn-r}
Let $f$ be a $t$-relaxed coloring of $K(n_1,n_2,\dots,n_s)$. If a color $i$ appears on exact $r\ge 2$ parts of  $K(n_1,n_2,\dots,n_s)$, then $|f^{-1}(i)|\le t+ \lfloor \frac{t}{r-1} \rfloor$.
\end{lemma}

\begin{pf}
Suppose the color $i$ appears on the $r$ parts $V_{j_1},V_{j_2},\dots,V_{j_r}$. For $q=1,2,\dots,r$, let $x_q=|V_{j_q}\cap f^{-1}(i)|$. Since $f$ is a $t$-relaxed coloring, for $p=1,2,\dots,r$, we have $\sum\limits_{q=1}^{r} x_q-x_p\le t$. Summing up all these $r$ inequalities together, we get $(r-1)\sum\limits_{q=1}^{r} x_q \le rt$. Thus $|f^{-1}(i)|=\sum\limits_{q=1}^{r} x_q \le t+ \frac{t}{r-1}$. Note that $|f^{-1}(i)|$ is an integer, the lemma follows.
\end{pf}

\begin{corollary} \label{beta}
Let $K(n_1,n_2,\dots,n_s)$ be the complete $s$-partite with $n_1\ge n_2\ge \cdots \ge n_s$. If $n\ge 2t$ then $\beta_t(K(n_1,n_2,\dots,n_s))=n_1$, and if $n_1\le 2t$ then $\beta_t(K(n_1,n_2,\dots,n_s))\le 2t$.
\end{corollary}

\begin{corollary} \label{geq2t}
If $n_j\ge 2t$ for each $j=1,2,\dots,s$, then $\chi_t(K(n_1,n_2,\dots,n_s))=s$.
\end{corollary}

\begin{pf}
It is clear that $\chi_t(K(n_1,n_2,\dots,n_s))\le s$.
We prove the opposite inequality by induction on $s$. It is obviously true for $s=1$. Suppose it holds for all complete $(s-1)$-partite graphs with each part having at least $2t$ vertices. Now consider the graph $K(n_1,n_2,\dots,n_s)$ with $n_s\ge 2t$. 

Let $f$ be any $t$-relaxed coloring of $K(n_1,n_2,\dots,n_s)$. If a color $i$ appears only on one part, say $V_j$,  then by the induction hypothesis, $\chi_t(K(n_1,n_2,\dots,n_s)-V_j)\ge s-1$, implying that $||f||\ge s$. Thus we assume that each color appears on at least $2$ parts. By Lemma \ref{mvn-r}, each color is assigned to at most $2t$ vertices of $K(n_1,n_2,\dots,n_s)$, implying that $||f||\ge s$. Therefore $\chi_t(K(n_1,n_2,\dots,n_s))\ge s$ and the corollary follows.
\end{pf}

This corollary provides a class of graphs $G$ with $\chi_t(G)=\chi(G)$ for any positive integer $t$. 

\begin{theorem} \label{2t}
Suppose $G=K(n_1,n_2,\dots,n_s)$ has $r$ parts of order at least $2t$. Let $\sigma_{r+1,s}=\sum\limits_{j=r+1}^{s} n_j$. Then $r+\lceil \frac{\sigma_{r+1,s}}{2t} \rceil \le \chi_t(G)\le r+\lceil \frac{\sigma_{r+1,s}}{t+1} \rceil$. Both the lower and upper bounds are attainable. 
\end{theorem}

\begin{pf}
We first prove the upper bound. Let $H_1=K(n_{1},\dots,n_r)$ and $H_2=K(n_{r+1},\dots,n_s)$.  It is obvious that $\chi_t(G)\le \chi_t(H_1)+\chi_t(H_2)$.  By Corollary \ref{geq2t}, $\chi_t(H_1)=r$. On the other hand, by Corollary \ref{ub}, $\chi_t(H_2)\le \lceil \frac{\Delta(H_2)+1}{t+1}\rceil \le \frac{|V(H_2)|}{t+1} = \lceil \frac{\sigma_{r+1,s}}{t+1} \rceil$. Thus $\chi_t(G)\le r+\lceil \frac{\sigma_{r+1,s}}{t+1} \rceil$.

We now show the lower bound. Suppose $f$ is an optimal $t$-relaxed coloring of $G$. Let
\[
I=\{i:\ \exists \ j\in \{1,2,\dots,r\} \ \mbox{s.t. $|f^{-1}(i)\cap V_j|\ge t+1$} \}.
\]
It is clear that, for each $i\in I$, there is an integer $j_i\in \{1,2,\dots,r\}$ such that $f^{-1}(i) \subseteq V_{j_i}$. Thus $\sum\limits_{i\in I} |f^{-1}(i)| \le  \sum\limits_{i\in I} |V_{j_i}| \le n_1+n_2+\cdots + n_{|I|}$. From the definition of $I$, if $i\in f(G)\setminus I$ then $|f^{-1}(i)\cap V_j| \le t$ for all $j\in \{1,2,\dots,r\}$, implying by Lemma \ref{mvn-r} that $|f^{-1}(i)| \le 2t$. Then
\[
\sum\limits_{i=1}^{s} n_i =\sum\limits_{i\in I} |f^{-1}(i)| + \sum\limits_{i\in f(G)\setminus I}  |f^{-1}(i)| 
\le n_1+n_2+\cdots +n_{|I|} + 2t |f(G)\setminus I|. 
\] 
It follows that $|f(G)\setminus I|\ge \frac{\sum\limits_{i=|I|+1}^{r} n_i }{2t} + \frac{\sum\limits_{i=r+1}^{s} n_i }{2t} \ge r-|I|+ \frac{\sum\limits_{i=r+1}^{s} n_i}{2t} $. Therefore $\chi_t(G)\ge r+\lceil \frac{\sigma_{r+1,s}}{2t} \rceil$.

Let $G_1=K(n_1,n_2,\dots,n_s)$ with $n_1=\cdots=n_r=2t+1$ and $n_{r+1}=\cdots=n_s=1$. Suppose $f$ is an optimal $t$-relaxed coloring of $G_1$. Let $I=f(V_1\cup\cdots\cup V_r)\setminus f(V_{r+1}\cup\cdots\cup V_s)$. It is clear that, for each $i\in I$, $|f^{-1}(i)|\le 2t+1$. For each $i\in f(G_1)\setminus I$, there exists some $j\in\{r+1,r+2,\dots,s\}$ such that $f^{-1}(i)\cap V_j\not=\emptyset$, implying that $|f^{-1}(i)|\le t+1$. Then $(2t+1)r+s-r\le (2t+1)|I| + (t+1) |f(G_1)\setminus I|$. On the other hand, in an optimal $t$-relaxed coloring of $G_1$, we must have $|I|\le r$.  Thus $|f(G_1)\setminus I| \ge \frac{(2t+1)(r-|I|)}{t+1} +\frac{s-r}{t+1} \ge r-|I| +\frac{s-r}{t+1}$, implying that $|f(G_1)|\ge r+\lceil \frac{\sigma_{r+1,s}}{t+1} \rceil$. Therefore  $\chi_t(G_1)= r+\lceil \frac{\sigma_{r+1,s}}{t+1} \rceil$, attaining the upper bound.

Let $G_2=K(n_1,n_2,\dots,n_s)$ with $n_1=\cdots=n_r=2t$ and $n_{r+1}=\cdots=n_s=t$.  Define a $t$-relaxed coloring $f$ of $G_2$ as follows. For each $j\in \{1,2,\dots,r\}$, let $f(V_j)=\{j\}$. For each $j\in \{1,2,\dots,\lceil \frac{s-r}{2} \rceil -1\}$, let $f(V_{r+2j-1}\cup V_{r+2j})=\{r+j\}$. Finally, if $s-r$ is even then let $f(V_{s-1}\cup V_{s})=\{\frac{s-r}{2}\}$, and if $s-r$ is odd then let $f(V_{s})=\{\lceil \frac{s-r}{2} \rceil\}$. It is clear that $|f(G_2)|=r+ \lceil \frac{s-r}{2} \rceil=r+ \lceil \frac{\sigma_{r+1,s}}{2t} \rceil$.  Thus $\chi_t(G_2)= r+\lceil \frac{\sigma_{r+1,s}}{2t} \rceil$, attaining the lower bound.
\end{pf}

\noindent{\bf Remark:}  Let $G=K(n_1,n_2,\dots,n_s)$ be a complete multi-partite graph with $s$ parts. Suppose $n_i<2t$ for each $i=1,2,\dots,s$. According to Theorem \ref{2t}, $\lceil \frac{n_1+n_2+\cdots+n_s}{2t} \rceil \le \chi_t(G)\le \lceil \frac{n_1+n_2+\cdots+n_s}{t+1} \rceil$. Let $n=n_1+n_2+\cdots+n_s$. Then, for each $k\in [\lceil \frac{n}{2t} \rceil, \lceil \frac{n}{t+1} \rceil]$, to determine if $G$ is $(k,t)$-colorable, the algorithm given by Belmonte etc. in [\citeonline{BLM2017}] runs in $O^* (k^{O(t^4)})$ time steps. Even when $t=3$, the time complexity is about $(\frac{n}{6})^{81}$.  It is thus meaningful to design more efficient algorithm for computing the $t$-relaxed chromatic numbers of complete multi-partite graphs.  We are trying to do so in the following sections.

\section{The greedy algorithm}

In this section, we present a greedy algorithm to construct a $t$-relaxed coloring of complete multi-partite graphs. Note that the deletion of any subset of vertices of a complete multi-partite graph results in a new complete multi-partite graph. The key idea of the greedy algorithm is that each time we find a maximum $t$-sparse set of the subgraph induced by all uncolored vertices and assign a new color to all vertices in it. We need to specify the technique how to find a maximum $t$-sparse set of the subgraph induced by all uncolored vertices so that the complexity of the greedy algorithm can be estimated accurately. 

In general, the problem of deciding if $\beta_t(G)\ge k$ is NP-complete, see [\citeonline{LY1980}]. However, in this section, we shall first show that $\beta(G)$ for multi-partite graphs can be computed in $O((t+1)^2)$ time steps, provided that the order of parts are given in non-increasing order.  

Let the $s$ parts of $G=K(n_1,n_2,\dots,n_s)$ be $V_1,V_2,\dots,V_s$. For $j\in\{1,2,\dots,s\}$, let $V_j=\{v_j^1,v_j^2,\dots,v_j^{n_j}\}$.  Sort the parts of $G$ such that $n_1\ge n_2\ge \cdots \ge n_s$. If $n_1\ge 2t$ then according to Corollary \ref{beta} the set $V_1$ is the largest $t$-sparse set of $G$. Thus we assume $n_1<2t$. 

Note that any $t$-sparse set of $G$ intersects with at most $t+1$ parts. To find a maximum $t$-sparse set of $G$, it suffices to find, for every $i=2,3,\dots,\min\{s,t+1\}$,  a maximum $t$-sparse set of $G$ with the property that it intersects exactly $i$ parts of $G$. The maximum is clearly attained by the first $i$ parts. Thus we only need to find a maximum $t$-sparse set of $G$ with the property that it intersects exactly the first $i$ parts of $G$. For $j=1,2,\dots,i$, let $x_j$ be the number of vertices we are going to choose from the part $V_j$. The objective is to maximize the sum $x_1+x_2+\cdots+x_i$ subject to $1\le x_j\le n_j$ and $\sum\limits_{h=1}^{i}x_h-x_j\le t$ for each $j=1,2,\dots,i$. Again, without loss of generality, we may assume $x_1\ge x_2\ge \cdots \ge x_i$.  For $i\in\{2,3,\dots,\min\{s,t+1\}\}$, we define this integral linear programming, denote by $LP(i)$, as follows:

\[
\begin{array}{ll}
LP(i):  & \max: \ x_1+x_2+\cdots+x_i\\
s.t.  &  \sum\limits_{h=1}^{i}x_h-x_j\le t, \ j=1,2,\dots,i;\\
       & 1\le x_j\le n_j, \ j=1,2,\dots,i;\\
       &  x_1\ge x_2\ge \cdots \ge x_i.
\end{array}
\]

If $n_i\le \lfloor \frac{t}{i-1} \rfloor$ then let $x_i=n_i$ and choose $(x_1,x_2,\dots,x_{i-1})$ 
with $x_1\ge x_2\ge \cdots \ge x_{i-1}\ge n_i$ so that $x_1+x_2\cdots+x_{i-1}=\min \{t,n_1+\cdots+n_{i-1}\}$. 
If $n_i > \lfloor \frac{t}{i-1} \rfloor$ then let $x_i=\lfloor \frac{t}{i-1}\rfloor$ and choose $(x_1,x_2,\dots,x_{i-1})$ with $x_1\ge x_2\ge \cdots \ge x_{i-1}\ge n_i$  so that $x_1+x_2\cdots+x_{i-1}=t$.  It is straightforward to check that this gives an optimal solution of $LP(i)$.\\

\begin{problem}\label{max-t-sparse}\ 
{\em (Find a maximum $t$-sparse set of a complete multi-partite graph)}\\
\noindent {\em Input}: A complete multi-partite graph $G$.\\
\vskip -0.8cm
\noindent {\em Output}: A maximum $t$-sparse set of $G$.
\end{problem}

According to above discussions, the algorithm to compute the maximum $t$-sparse set of a complete multi-partite graph can be describe as in Algorithm 1.
 
\begin{algorithm}[h!]\label{max-t-sparse set}
\caption{Find a maximum $t$-sparse set of a complete multi-partite graph. }
\begin{algorithmic}[1]
\STATE Sort the parts of $G$ such that $n_1\ge n_2\ge \cdots \ge n_s$.\\
\STATE Set $s^*:=\min\{t+1,s\}$.
\STATE {\bf If} $n_1\ge 2t$ {\bf then} output $x^{(1)}=(n_1)$.
\STATE {\bf else} for $i=2$ to $s^*$ {\bf do} 
\STATE \hskip 1cm Solve $LP(i)$. \\
             \hskip 1cm Set $x^{(i)}$ be the optimal solution.\\ 
             \hskip 1cm Set $M_i$ be the optimal value of $LP(i)$. 
\STATE Find the minimum number $k$ such that $M_k=\max\{M_2,M_3,\dots,M_{s^*}\}$.
\STATE Output $x^{(k)}$.
\end{algorithmic}\label{max-t-sparse set}
\end{algorithm}

We are now ready to present the greedy algorithm for $t$-relaxed coloring of a complete multi-partite graph.\\

\begin{problem}\label{opt-t-relaxed-col}\
{\em (Find an optimal $t$-relaxed coloring of a complete multi-partite graph)}\\
\noindent {\em Input}: A complete multi-partite graph $G$.\\
\vskip -0.8cm
\noindent {\em Output}: An optimal $t$-relaxed coloring of $G$.
\end{problem}

For a subset $S$ of $V(G)$, we use $G[S]$ to denote the subgraph of $G$ induced by $S$.\\

\begin{algorithm}[h!]\label{t-re-c}
\caption{Construct a $t$-relaxed coloring of a complete multi-partite graph}
\begin{algorithmic}[1]
\STATE Set $V=V_1\cup V_2\cup \cdots \cup V_s$.
\STATE Set $i:=1$ and $R:=\emptyset$.\\
\STATE {\bf While}  $V\setminus R$ is not empty {\bf do}:  \\
             \hskip 1cm Use Algorithm \ref{max-t-sparse set} to find a maximum $t$-sparse set $W$ of $G[V\setminus R]$. \\
             \hskip 1cm Assign the color $i$ to all vertices in $W$.\\
             \hskip 1cm Set $R:=R\cup W$ and $i:=i+1$.
\STATE Output the coloring.
\end{algorithmic}\label{t-re-c}
\end{algorithm}

It is easy to see that Algorithm \ref{max-t-sparse set} runs in $O((t+1)^2)$ time. While Algorithm \ref{t-re-c} will call Algorithm \ref{max-t-sparse set} at most $\lceil \frac{n}{t+1} \rceil$ times. Thus the overall complexity of Algorithm \ref{t-re-c} is $O(tn)$.

Suppose $G=K(n_1,n_2,\dots,n_s)$ has $r$ parts of order at least $2t$. Let $\sigma_{r+1,s}=\sum\limits_{j=r+1}^{s} n_j$. Then by  Theorem \ref{2t}, we have $r+\lceil \frac{\sigma_{r+1,s}}{2t} \rceil \le \chi_t(G)\le r+\lceil \frac{\sigma_{r+1,s}}{t+1} \rceil$. According to the steps of  Algorithm \ref{t-re-c}, the $t$-relaxed coloring produced by Algorithm \ref{t-re-c} will use at most $r+\lceil \frac{\sigma_{r+1,s}}{t+1} \rceil$ colors. Thus the approximation ratio of Algorithm \ref{t-re-c} is at most $\frac{r+\lceil \frac{\sigma_{r+1,s}}{t+1} \rceil}{r+\lceil \frac{\sigma_{r+1,s}}{2t} \rceil}\le 2$. 

In what follows, we first prove that for $t\in \{1,2,3,4,5,6\}$, the greedy algorithm produces an optimal $t$-relaxed coloring of a complete multi-partite graph. And then, for $t\ge 7$, we construct examples  to illustrate that the greedy strategy does not always construct an optimal $t$-relaxed coloring. We only present the proofs of $t=1$ and $t=3$. Though the idea behind the proofs are the same, the proofs of $t=4,5,6$ are somewhat complicated than that of $t=3$. So we omit them here.

\section{Determining the $1$-relaxed chromatic number of $K(n_1,n_2,\dots,n_s)$}

\begin{lemma} \label{t=1}
Let $j$ be an integer in $\{1,2,\dots,s\}$. If $n_j\ge 2$ then there is an optimal $1$-relaxed coloring $f$ of $K(n_1,n_2,\dots,n_s)$ with $|f(V_j)|=1$.
\end{lemma}

\begin{pf}
Without loss of generality, we may assume that $n_j=2$. Let $f$ be an optimal $1$-relaxed coloring of $K(n_1,n_2,\dots,n_s)$. If $|f(V_j)|=1$ then we are done.

Suppose $|f(V_j)|=2$. Let $v_j^1$ and $v_j^2$ be the two vertices in $V_j$ with $f(v_j^1)=i_1$ and  $f(v_j^2)=i_2$, where $i_1\not=i_2$. If the color $i_1$ (resp. $i_2$) does not appear on any vertex outside $V_j$, then, by recoloring the vertex $v_j^2$ (resp. $v_j^1$)  with the color $i_1$ (resp. $i_2$), we get a $1$-relaxed coloring of $K(n_1,n_2,\dots,n_s)$. Denote this coloring by $f'$. It is clear that $f'$ is optimal and $|f'(V_j)|=1$.

Thus we assume that there are two vertices $w_1$ and $w_2$ outside $V_j$ with $f(w_1)=i_1$ and $f(w_2)=i_2$. Since $f$ is a $1$-relaxed coloring, $f^{-1}(i_1)=\{v_j^1,w_1\}$ and $f^{-1}(i_2)=\{v_j^2,w_2\}$. Now, by recoloring the vertex $w_1$ with the color $i_2$ and the vertex $v_j^2$ with the color $i_1$, one gets an optimal $1$-relaxed coloring of $K(n_1,n_2,\dots,n_s)$. Denote the so obtained coloring by $f''$. It is clear that $f''$ is optimal and $|f''(V_j)|=1$. The lemma follows.
\end{pf}

We immediately have the following corollary.

\begin{corollary} \label{+1}
Let $j$ be an integer in $\{1,2,\dots,s\}$. If $n_j\ge 2$ then $\chi_1(K(n_1,n_2,\dots,n_s))=\chi_1(K(n_1,n_2,\dots,n_s)-V_j)+1$.
\end{corollary}

Since $\chi_t(K_n)=\lceil \frac{n}{t+1} \rceil$, Corollary \ref{+1} implies the following theorem.

\begin{theorem}  \label{1-relax}
Let $G$ be any complete $s$-partite graph. Suppose there are $r$ single vertex parts in $G$. Then $\chi_1(G)=s-r+\lceil \frac{r}{2} \rceil$.
\end{theorem}

It follows from the theorem that Algorithm \ref{t-re-c} presented in the previous section will output an optimal $1$-relaxed coloring of a complete multi-partite graph.

\section{Computing the $3$-relaxed chromatic number of $K(n_1,n_2,\dots,n_s)$}

In this section, we prove that the greedy algorithm produces an optimal $t$-relaxed coloring of a complete multi-partite graph for $t\le 6$. This would be proved one by one for $t=2,3,4,5,6$. However the idea behind the proof is essentially the same. Thus, in the following, we only present the proof for $t=3$.

Let $V=\cup_{j=1}^{s} V_j$ be the vertex set of $K(n_1,n_2,\dots,n_s)$. For an integer $j\in \{1,2,\dots,s\}$, the $n_j$ vertices in part $V_j$ are denoted by $v_j^1,v_j^2,\dots,v_j^{n_j}$.
Suppose $f$ is a $t$-relaxed coloring of $K(n_1,n_2,\dots,n_s)$. Then $f(V)$ is the set of colors used by $f$. For $r\in f(V)$ and $h\in \{1,2,\dots,s\}$, let $V_{r,h}$ denote the set of vertices in $V_h$ that are assigned the color $r$  under $f$. And let  $m_{r,h}$ be the order of $V_{r,h}$. 

\begin{lemma} \label{3-2t}
Suppose $G=K(n_1,n_2,\dots,n_s)$ has at least one part having at least $6$ vertices. Then there is an optimal $3$-relaxed coloring $f$ and a part $V_j$ with $n_j\ge 6$ such that $|f(V_j)|=1$. 
\end{lemma}

\begin{pf}
Let $f$ be an optimal $3$-relaxed coloring of $K(n_1,n_2,\dots,n_s)$. Choose $i\in f(V)$ and $j\in \{1,2,\dots,s\}$ such that $m_{i,j}=\max\{m_{r,h}:\  r\in f(V)\ \mbox{and $h\in\{1,2,\dots,s\}$ with $n_h\ge 6$}\}$. We consider the following four cases.

Case 1. $m_{i,j}\ge 4$.
Since $f$ is a $3$-relaxed coloring of $G$, in this case, we have $m_{i,p}=0$ for each $p\in \{1,2,\dots,s\}\setminus \{j\}$. If there are some vertices in $V_j$ that are not colored by the color $i$, then we can recolor them by the color $i$. And the resulting coloring is the desired one.

Case 2. $m_{i,j}=3$. If $f^{-1}(i)\cap (V\setminus V_j)=\emptyset$ then we can deal with it as in Case 1. Otherwise, there is exactly one other part of $G$, say $V_p$, with $V_{i,p}\not=\emptyset$. Let $V_{i,p}=\{v_p^1,\dots,v_p^b\}$ for some $b\in \{1,2,3\}$. Then we define a new $3$-relaxed coloring $f'$ as follows
\[
\begin{array}{ll}
f'(v_j^d)=i, ~\mbox{for}~ 1\le d \le n_j, \\
f'(v_p^d)=f(v_j^{3+d}), ~\mbox{for}~ 1\le d\le b, \\
f'(u)=f(u), ~\mbox{for all other vertices  $u$ of}~ G.
\end{array}
\]
It is easy to check that $f'$ is an optimal $3$-relaxed coloring of $G$ with $|f'(V_j)|=1$. 

Case 3. $m_{i,j}=2$. Then $m_{i^*,j}\le 2$ for any $i^*\in f(V)$. In this case, there are at most two other parts of $G$ with vertices assigned the color $i$. If there is at most one other part of $G$ with vertices assigned the color $i$, then we can deal with it as in Case 2. Thus we assume that there are exact two other such parts of $G$, say $V_p$ and $V_q$. Since $f$ is a $3$-relaxed coloring, we must have $m_{i,p}=m_{i,q}=1$. Let $v_p^1$ and $v_q^1$ be the two vertices in $V_{i,p}$ and $V_{i,q}$. Since $m_{i,j}=2$ and $n_j\ge 6$, except the color $i$, there are at least two other distinct colors, say $i_1$ and $i_2$, in $f(V_j)$. We now define a new $3$-relaxed coloring $f'$ as follows
\[
\begin{array}{ll}
f'(v_j^d)=i, ~\mbox{for}~ 1\le d \le n_j, \\
f'(v_p^1)=i_1, ~\mbox{and}~ f'(v_q^1)=i_2,\\
f'(u)=f(u), ~\mbox{for all other vertices  $u$ of}~ G.
\end{array}
\]
It is easy to check that $f'$ is an optimal $3$-relaxed coloring of $G$ with $|f'(V_j)|=1$. 

Case 4. $m_{i,j}=1$. Then $m_{i^*,j}\le 1$ for any $i^*\in f(V)$. Let $b$ be the number of vertices in $V\setminus V_j$ that are assigned the color $i$ and let $w_1,w_2,\dots,w_b$ be these $b$ vertices.  Since $f$ is a $3$-relaxed coloring of $G$,  we have $b\le 3$. Note that $n_j\ge 6$ and $m_{i^*,j}\le 1$ for any $i^*\in f(V)$, there are at least $5$ different colors, say $i_1,i_2,\dots,i_5$, each appearing on exact one vertex in $V_j$. We now define a new $3$-relaxed coloring $f'$ as follows
\[
\begin{array}{ll}
f'(v_j^d)=i, ~\mbox{for}~ 1\le d \le n_j, \\
f'(w_d)=i_d, ~\mbox{for}~ 1\le d \le b,\\
f'(u)=f(u), ~\mbox{for all other vertices  $u$ of}~ G.
\end{array}
\]
It is easy to check that $f'$ is an optimal $3$-relaxed coloring of $G$ with $|f'(V_j)|=1$. 
\end{pf}

By Lemma \ref{3-2t}, we immediately have the following corollary.

\begin{corollary} \label{3-2t+1}
Let $j$ be an integer in $\{1,2,\dots,s\}$. If $n_j\ge 6$ then $\chi_3(K(n_1,n_2,\dots,n_s))=\chi_3(K(n_1,n_2,\dots,n_s)-V_j)+1$.
\end{corollary}

From Corollary \ref{3-2t+1}, to compute the $3$-relaxed chromatic number of a complete multi-partite graph it suffices to consider those complete multi-partite graphs with each part having at most $5$ vertices. Thus, from now on, we assume that $G$ is a complete multi-partite graph $K(n_1,n_2,\dots,n_s)$ with $n_j\le 5$ for $j=1,2,\dots,s$.

Let $f$ be an optimal $3$-relaxed coloring of $K(n_1,n_2,\dots,n_s)$. For convenience, for each part $V_j$, we relabel all vertices in $V_j$ according to the colors assigned to them. After the relabeling, for any two vertices $v_j^x$ and $v_j^y$ in $V_j$ with $x<y$, if $f(v_j^x)=f(v_j^y)$ then $f(v_j^q)=f(v_j^x)$ for all $q\in \{x,x+1,\dots,y\}$, and if $f(v_j^x)\not=f(v_j^y)$ then $|f^{-1}(f(v_j^x))|\ge |f^{-1}(f(v_j^y))|$. For example, suppose $f(V_j)=\{1,2,4\}$ and $|f^{-1}(1)|=2$, $|f^{-1}(2)|=1$ and $|f^{-1}(4)|=3$, then $f^{-1}(4)=\{v_j^1,v_j^2,v_j^3\}$, $f^{-1}(1)=\{v_j^4,v_j^5\}$ and $f^{-1}(2)=\{v_j^6\}$.

\begin{lemma} \label{3-33}
Let $G$ be a complete multi-partite graph $K(n_1,n_2,\dots,n_s)$ with $n_j\le 5$ for $j=1,2,\dots,s$.
Suppose $G$ has at least two parts of order at least $3$. Then there is an optimal $3$-relaxed coloring of $G$  which assigns three vertices of $V_1$ and three vertices of $V_2$ the same color. 
\end{lemma}

\begin{pf}
Let $f$ be an optimal $3$-relaxed coloring of $K(n_1,n_2,\dots,n_s)$. Choose $i\in f(V)$ and $j\in \{1,2\}$ such that $m_{i,j}=\max\{m_{r,h}: r\in f(V)\ \mbox{and $h\in\{1,2\}$}\}$. Let $\{j,q\}=\{1,2\}$. We consider the following four cases.

Case 1. $m_{i,j}\ge 4$. Then $m_{i,p}=0$ for each $p\in \{q\}\cup \{3,4,\dots,s\}$. We now define a new coloring $f'$ as follows
\[
\begin{array}{ll}
f'(v_j^d)=i, ~\mbox{for}~ 1\le d \le 3, \\
f'(v_q^d)=i, ~\mbox{for}~ 1\le d \le 3, \\
f'(v_j^d)=f(v_q^{d-3}), ~\mbox{for}~ 4\le d \le n_j,\\
f'(u)=f(u), ~\mbox{for all other vertices $u$ of}~ G.
\end{array}
\]
It is clear that $|f'(V)|\le |f(V)|$. If $f(v_q^1)\not=f(v_q^4)$ then we also have $f(v_q^2)\not=f(v_q^4)$. It is easy to see that $f'$ is the desired $3$-relaxed coloring of $G$. If $f(v_q^1)=f(v_q^4)$ then  $f(v_q^2)=f(v_q^3)=f(v_q^4)$. It follows that the color $f(v_q^1)$ only appears on vertices in $V_q$, implying that $f'$ is the desired $3$-relaxed coloring of $G$.

Case 2. $m_{i,j}=3$. Then $m_{i,q}\le 3$ and there is at most one integer $p\in \{q\}\cup \{3,4,\dots,s\}$ with $m_{i,p}>0$. If  $m_{i,p}=0$ for all $p\in \{q\}\cup \{3,4,\dots,s\}$ then by recoloring the three vertices $v_q^1,v_q^2,v_q^3$ with the color $i$ we obtain a desired $3$-relaxed coloring of $G$. Suppose $m_{i,q}>0$. If $m_{i,q}=3$ then we are done. If $m_{i,q}<3$ then choose any $3-m_{i,q}$ vertices in $V_q$ that are not colored by $i$ and recolor them with the color $i$, resulting a desired $3$-relaxed coloring of $G$.

Now suppose $m_{i,q}=0$ and there is an integer $p\in \{3,4,\dots,s\}$ with $m_{i,p}>0$. Then $1\le m_{i,p}\le 3$. Let $m_{i^*,q}=\max\{m_{r,q}: r\in f(V_q)\}$. If $m_{i^*,q}=3$ then recolor the three vertices $v_q^1,v_q^2,v_q^3$ with the color $i$ and recolor the vertices in $V_{i,p}$ with the color $i^*$. It is not difficult to check that the resulting coloring is the desired one. If $m_{i^*,q}=1$ then recolor the three vertices $v_q^1,v_q^2,v_q^3$ with the color $i$ and recolor the vertices in $V_{i,p}$ with different colors from $\{ f(v_q^1),f(v_q^2),f(v_q^3) \}$. Again we can check that the resulting coloring is the desired one. 

We are now left with the case $m_{i^*,q}=2$. Let $V_{i^*,q}=\{v_q^1,v_q^2\}$. If $m_{i,p} \le 2$ then recolor the three vertices $v_q^1,v_q^2,v_q^3$ with the color $i$ and recolor the vertices in $V_{i,p}$ with the color $i^*$. It is straightforward to check that the resulting coloring is the desired one. Thus we assume $m_{i,p}=3$. If there is a color $i'$ appearing on exact one vertex, say $v_q^d$, in $V_q$, then recolor the three vertices $v_q^1,v_q^2,v_q^3$ with the color $i$ and recolor the three vertices in $V_{i,p}$ with the colors $i^*$ and $i'$, two of them with $i^*$ and the other with $i'$. It can be checked that the resulting coloring is the desired one. The remaining case now is that each color appears on exact two vertices in $V_q$, implying $n_q=4$. Let $i^*$ and $i'$ be the two colors appearing on vertices in $V_q$ with $V_{i^*,q}=\{v_q^1,v_q^2\}$ and $V_{i',q}=\{v_q^3,v_q^4\}$. If $|f^{-1}(i')\setminus V_q| \le 2$ then recolor the three vertices $v_q^1,v_q^2,v_q^3$ with the color $i$ and recolor the three vertices in $V_{i,p}$ with the colors $i^*$ and $i'$, two of them with $i^*$ and the other with $i'$. It is easy to see that the resulting coloring is the desired one. If $|f^{-1}(i')\setminus V_q| =3$ then there is a part $V_l$ with $f^{-1}(i')\setminus V_q=V_{i',l}$. Recolor the three vertices $v_q^2,v_q^3,v_q^4$ with the color $i$ and recolor the three vertices in $V_{i,p}$ with the same color $i'$. It is not difficult to check that the resulting coloring is the desired one. 

Case 3. $m_{i,j}=2$. Let $V_{i,j} =\{v_j^1,v_j^2\}$. 
Then $m_{r,h}\le 2$ for all $r\in f(G)$ and $h\in \{1,2\}$ and there are at most two integers $p\in \{q\} \cup \{3,4,\dots,s\}$ with $m_{i,p}>0$. Furthermore, if $p$ and $p'$ are two integers in $\{q\} \cup \{3,4,\dots,s\}$ with $m_{i,p}>0$ and $m_{i,p'}>0$ then $m_{i,p}=m_{i,p'}=1$. It is clear that $|f^{-1}(i)\setminus V_j|\le 3$ and if $|f^{-1}(i)\setminus V_j|= 3$ then there is some $p\in \{1,2,\dots,s\}\setminus \{j\}$ with $V_{i,p}=f^{-1}(i)\setminus V_j$.

If $f(v_j^3)\not=f(v_j^4)$ then recolor the vertex $v_j^3$ with the color $i$ and recolor a vertex in 
$f^{-1}(i)\setminus V_j$ with the color $f(v_j^3)$. It is easy to see that the resulting coloring is also an optimal $3$-relaxed coloring of $G$ and we are back to Case 2. 

Now suppose $f(v_j^3)=f(v_j^4)=i^*$. If $|f^{-1}(i^*)\setminus V_j|\le 2$ then we can do the same thing as in the previous paragraph and go back to Case 2. Thus we assume $|f^{-1}(i^*)\setminus V_j|=3$ and there is a part $V_p$ with $V_{i^*,p}=f^{-1}(i^*)\setminus V_j$. If $|f^{-1}(i)\setminus V_j|\le 2$ then we let the color $i^*$ replace the role of the color $i$ and deal with it as in the previous paragraph. Therefore we assume $|f^{-1}(i)\setminus V_j|=3$ and there is a part $V_q$ with $V_{i,q}=f^{-1}(i)\setminus V_j$. Then recolor all vertices in $V_{i^*,j}$ with the color $i$ and all vertices in $V_{i,q}$ with the color $i^*$ and go back to Case 1.

Case 4. $m_{i,j}=1$. Then $m_{i^*,h}\le 1$ for any $i^*\in f(V)$ and $h\in \{j,q\}$. Let $f(v_j^1)=i$ and $f(v_j^2)=i^*$. It is obvious that $|f^{-1}(i^*)\setminus V_j|\le 3$. Then, by recoloring the vertex $v_j^2$ with the color $i$ and a vertex in $f^{-1}(i)\setminus V_j$ with the color $i^*$, we obtain a new optimal $3$-relaxed coloring of $G$ with $m_{i,j}=2$ and we are back to Case 3.
\end{pf}

By repeatedly applying Lemmas \ref{3-2t} and \ref{3-33}, one can reduce the problem to construct an optimal $3$-relaxed coloring of $K(n_1,n_2,\dots,n_s)$ to the case when $n_j\le 2$ for all $j\in \{1,2,\dots,s\}$ or the case when $3\le n_1\le 5$ and $n_j\le 2$ for all $j\in \{2,3,\dots,s\}$. These cases are going to be discussed in the following three lemmas. 

\begin{lemma} \label{3-32}
Let $G$ be a complete multi-partite graph $K(n_1,n_2,\dots,n_s)$ with $3\le n_1\le 5$ and $n_2=2$. Then there is an optimal $3$-relaxed coloring of $G$ which assigns three vertices of $V_1$ and the two vertices of $V_2$ the same color. 
\end{lemma}

\begin{pf}
Let $f$ be an optimal $3$-relaxed coloring of $K(n_1,n_2,\dots,n_s)$. Choose the color $i$ such that $m_{i,1}=\max\{m_{r,1}:\  r\in f(V_1)\}$.  We consider the following four cases.

Case 1. $m_{i,1}\ge 4$. Then $m_{i,p}=0$ for all $p\in\{2,,3,\dots,s\}$. We can recolor the two vertices in $V_2$ with the color $i$ and recolor the vertices $v_1^d$ for $4\le d\le m_{i,1}$ with colors from $\{f(v_2^1),f(v_2^2)\}$. It is easy to see that the resulting coloring is the desired one.

Case 2. $m_{i,1}=3$. Then there is at most one $p\in \{2,3,\dots,s\}$ with $m_{i,p}>0$. If $m_{i,p}=0$ for $p\in \{3,4,\dots,s\}$ then we are done by just recoloring all vertices in $V_2\setminus V_{i,2}$ with the color $i$. If there is some $p\in \{3,4,\dots,s\}$ with $m_{i,p}>0$ then recolor the two vertices in $V_2$ with the color $i$ and recolor the vertices in $V_{i,p}$ with colors from $\{f(v_2^1),f(v_2^2)\}$. It is clear that the resulting coloring is what we are desiring. 

Case 3. $m_{i,1}=2$.  Since $n_p\le 2$ for all $p\in \{2,3,\dots,s\}$, we must have $|f^{-1}(i)\setminus V_1|\le 2$. Then recolor the vertex $v_1^3$ with the color $i$ and recolor a vertex in $f^{-1}(i)\setminus V_1$ with the color $f(v_1^3)$. It is not difficult to check that the resulting coloring is also an optimal $3$-relaxed coloring of $G$ and we are back to Case 2. 

Case 4. $m_{i,1}=1$. Let $V_{i,1}=\{v_1^1\}$. It is obvious that $|f^{-1}(i)\setminus V_1|\le 3$. By recoloring the vertex $v_1^2$ with the color $i$ and  a vertex in $f^{-1}(i)\setminus V_1$ with the color $f(v_1^2)$, we obtain a new optimal $3$-relaxed coloring of $G$ with $m_{i,1}=2$ and we are back to Case 3. 
\end{pf}


\begin{lemma} \label{3-54-31}
Let $G$ be a complete multi-partite graph $K(n_1,n_2,\dots,n_s)$ with $3\le n_1\le 5$ and $n_2=1$. 
If $n_1\ge 4$ then there is an optimal $3$-relaxed coloring $f$ of $G$ with $|f(V_1)|=1$. And if $n_1=3$ then there is an optimal $3$-relaxed coloring of $G$ which assigns three vertices of $V_1$ and the only vertex of $V_2$ the same color. 
\end{lemma}

\begin{pf}
Let $f$ be an optimal $3$-relaxed coloring of $K(n_1,n_2,\dots,n_s)$. Choose the color $i$ such that $m_{i,1}=\max\{m_{r,1}:\  r\in f(V_1)\}$.  Let $V_{i,1}=\{v_1^1,\dots,v_1^{m_{i,1}}\}$.

Suppose $4\le n_1\le 5$. If $m_{i,p}=0$ for all $p\in \{2,3,\dots,s\}$ then we can color all vertices in $V_1$ with the color $i$. Thus we assume $1\le m_{i,1}\le 3$ and $f^{-1}(i)\setminus V_1\not=\emptyset$. Then it is clear that $h=|f^{-1}(i)\setminus V_1|\le 4-m_{i,1}$. Let the $h$ vertices in $f^{-1}(i)\setminus V_1$ be $w_1,w_2,\dots,w_h$. By recoloring the vertices $v_1^d$ ($m_{i,1}+1\le d\le n_1$) with the color $i$ and $w_g$ ($1\le g\le h$) with the color $f(v_1^{m_{i,1}+g})$, we obtain an optimal $3$-relaxed color of $G$ that assigns all vertices in $V_1$ the same color.

Suppose $n_1=3$.  If $1\le m_{i,1}\le 2$ then similar as above we can recolor all vertices in $V_1$ with the same color $i$ and left at most one vertex with color $i$ not in $V_1$. And we can easily obtain the desired coloring.
\end{pf}


\begin{lemma} \label{3-54-31}
Let $G$ be a complete multi-partite graph $K(n_1,n_2,\dots,n_s)$ with $1\le n_j\le 2$ for all $j\in\{1,2,\dots,s\}$. Then $\chi_3(G)=\lceil \frac{n_1+n_2+\dots+n_s}{4}\rceil$.   
\end{lemma}

\begin{pf}
Let $f$ be an optimal $3$-relaxed coloring of $G$. By Lemma \ref{mvn-r}, since $1\le n_j\le 2$  for all $j\in\{1,2,\dots,s\}$, it is clear that $f$ assigns each color to at most $4$ vertices of $G$. Thus $\chi_3(G)\ge \lceil \frac{n_1+n_2+\dots+n_s}{4}\rceil$. On the other hand, since $\beta_t(G)\ge \min\{t+1,|V(G)|\}$ for any graph $G$, the greedy algorithm uses at most $\lceil \frac{n_1+n_2+\dots+n_s}{4}\rceil$ colors.
\end{pf}


Let $G$ be a complete multi-partite graph $K(n_1,n_2,\dots,n_s)$ with $n_1\ge n_2\ge \cdots n_s$. Note that if $n_1\ge 6$ then $\beta_3(G)=6$, if $3\le n_2\le n_1\le 5$ then $\beta_3(G)=6$ and three vertices of $V_1$ together with three vertices of $V_2$ forms a maximum $3$-sparse set of $G$, if $3\le n_1 \le 5$ and $n_2=2$ then $\beta_3(G)=5$ and three vertices of $V_1$ together with the two vertices in $V_2$ forms a maximum $3$-sparse set of $G$, if $n_1\ge 4$ and $n_2=1$ then $\beta_3(G)=n_1$ and $V_1$ is a maximum $3$-sparse set of $G$, if $n_1 =3$ and $n_2=1$ then $\beta_3(G)=4$ and the three vertices of $V_1$ together with the vertex in $V_2$ forms a maximum $3$-sparse set of $G$, while if $n_1\le 2$ then $\beta_3(G)= \min\{4,|V(G)|\}$. Keeping in mind these facts, Lemmas and Corollaries in this section imply that Algorithm 2 always outputs an optimal $3$-relaxed coloring of a complete multi-partite graph. 

In next section, we shall illustrate by examples  that the same idea is not applicable to the case $t\ge 7$.

\section{Counterexamples for $t\ge 7$ to the greedy algorithm}

In this section, we present examples for which the greedy algorithm does not produce an optimal $t$-relaxed coloring of a complete multi-partite graph. We first deal with the case that $t$ is an odd number. 

\begin{example}\label{todd}
Let $t\ge 7$ be an odd number. Let $G$ be the graph $K(n_1,\dots,n_6)$ with $n_1=2t-1$, $n_2=t+1$ and $n_3=n_4=n_5=n_6=(t-1)/2$. Then $\chi_t(G)=3$. While the greedy algorithm uses $4$ colors.
\end{example}

From the structure of $G$, it is easy to observe that each color can be assigned to at most $2t$ vertices of $G$. Since $G$ has $5t-2>4t$ vertices, $\chi_t(G)\ge 3$. 

A $t$-relaxed coloring using $3$ colors is presented as follows. Use $g$ denote this coloring. $g^{-1}(i)$ is the set of vertices with color $i$. We use the notation $g^{-1}(i)=(x_1,x_2,x_3,x_4,x_5,x_6)$ to indicate that for $j\in\{1,2,\dots,6\}$ there are $x_j$ vertices in $V_j$ that are assigned the color $i$. Then the $t$-relaxed coloring $g$ using $3$ colors is defined as
\[
\begin{array}{llllllllllllll}
g^{-1}(1)=(2t-1,& 0,& 0,& 0,& 0,& 0& ),\\
g^{-1}(2)=(0,& \frac{t+1}{2},& \frac{t-1}{2},& \frac{t-1}{2},& 0,& 0&),\\
g^{-1}(3)=(0,& \frac{t+1}{2},& 0,& 0,& \frac{t-1}{2},& \frac{t-1}{2}& ).
\end{array}
\]
Thus $\chi_t(G)=3$.

Use $(n_1,n_2,n_3,n_4,n_5,n_6)$ to denote the complete multi-partite graph with the $j$th part having $n_j$ vertices for $j\in \{1,2,\dots,6\}$. Let $f$ be the coloring produced by the greedy algorithm. The algorithm works as follows.
\[
\begin{array}{llllllllllllll}
& & & & & & & (2t-1,& t+1,& \frac{t-1}{2},& \frac{t-1}{2},& \frac{t-1}{2},&\frac{t-1}{2}&  )\\
f^{-1}(1)=(t,& t,& 0,& 0,& 0,& 0& )  \rightarrow  & (t-1,& 1,& \frac{t-1}{2},& \frac{t-1}{2},& \frac{t-1}{2},& \frac{t-1}{2}&)\\
f^{-1}(2)=(\frac{t+1}{2},& 0,& \frac{t-1}{2},& \frac{t-1}{2},& 0,& 0&)  \rightarrow  &  (\frac{t-3}{2},&  1,&  0,&  0,&  \frac{t-1}{2},&  \frac{t-1}{2}&  )\\
f^{-1}(3)=(\frac{t-3}{2},& 0,& 0,& 0,& \frac{t-1}{2},& \frac{t-1}{2}& )   \rightarrow  & (0,&  1,&  0,&  0,&  0,&  0&  )\\
f^{-1}(4)=(0,& 1,& 0,& 0,& 0,& 0& )  \rightarrow  & (0,&  0,&  0,&  0,&  0,&  0&  )
\end{array}
\]
Thus the greedy algorithm uses $4$ colors.

\begin{example} \label{teven}
Let $t\ge 7$ be an even number. Let $G$ be the graph $K(n_1,\dots,n_6)$ with $n_1=2t-1$, $n_2=t+1$ and $n_3=n_4=t/2$, and $n_5=n_6=t/2-1$. Then $\chi_t(G)=3$. While the greedy algorithm uses $4$ colors.
\end{example}

The analysis of Example \ref{teven} is similar to that of Example \ref{todd}. So we omit the details.



\end{document}